\documentclass[a4paper,11pt]{amsart}
\usepackage[plainpages=false]{hyperref}
\usepackage{amsfonts,latexsym,rawfonts,amsmath,amssymb,amsthm, mathrsfs, lscape}
\usepackage{verbatim}

\usepackage[all]{xy}
\usepackage{graphicx,psfrag}

\usepackage{array, tabularx}

\usepackage{setspace}

\newtheorem{thm}{Theorem}[section]
\newtheorem{cor}[thm]{Corollary}
\newtheorem{lem}[thm]{Lemma}

\newtheorem{rmk}[thm]{Remark}

\theoremstyle{definition}
\newtheorem{defi}[thm]{Definition}

\numberwithin{equation}{section}

\def \C {\mathbb C}

\def \R {\mathbb R}
\def \D {\mathcal D}
\def \P {\mathbb P}
\def \X {\mathcal X}
\def \L {\mathcal L}
\def \p {\partial}
\def \bp {\bar{\partial}}
\def \bz {\bar{z}}
\def \bt {\bar{t}}
\def \bw {\bar{w}}

\begin{document}

\title{Note on K-stability of pairs}
\author{Song Sun}
\date{\today}
\address{Department of Mathematics, Imperial College, London SW7 2AZ, United Kingdom.}
\email{s.sun@imperial.ac.uk}
\maketitle

\begin{abstract}
We prove that a pair $(X, D)$ with $X$ Fano and $D$ an anti-canonical divisor is K-unstable for negative angles, and is K-semistable for zero angle. 
\end{abstract}

\section{Introduction}

Let $X$ be a Fano manifold. It was first proposed by Yau \cite{Ya2} that finding K\"ahler-Einstein metrics on $X$ should be related to a certain algebro-geometric stability. In \cite{Ti}, the notion of K-stability was introduced by Tian. This has been conjectured to be equivalent to the existence of a K\"ahler-Einstein metric. One direction is essentially known, in a wider context of constant scalar curvature K\"ahler metrics \cite{Do1}. Namely, it is proved by Donaldson \cite{Do2} that the existence of a constant scalar curvature metric implies K-semistability. This was later strengthened by Stoppa \cite{St} to K-stability in the absence of continuous automorphism group, and by Mabuchi \cite{Ma} to K-polystability in general.\\

Recently  in \cite{Do4}(see also, \cite{Sz}, \cite{Li})  K-stability has been defined for a pair $(X, D)$, where $X$ is a Fano manifold and $D$ is a smooth anti-canonical divisor. The definition involves a parameter $\beta\in\R$. At least when $\beta\in(0,1]$, the K-stability of a pair $(X, D)$ with parameter $\beta$ is conjectured to be equivalent to the existence of a K\"ahler-Einstein metric on $X$ with cone singularities of angle $2\pi\beta$ transverse to $D$. This generalization grew out of a new continuity method for dealing with the other direction of the above conjecture, as outlined in \cite{Do3}.  Note   heuristically the case $\beta=0$ corresponds to a complete Ricci flat metric on the complement $X\setminus D$. By the work of Tian-Yau \cite{TY} such a metric always exists if $D$ is smooth. In this short article we prove the following theorem, which may be viewed as an algebraic counterpart of the differential geometric result of Tian-Yau. 

\begin{thm}\label{thm}
Any pair $(X, D)$ is strictly K-semistable with respect to angle $\beta=0$, and  K-unstable with respect to angle $\beta<0$.
\end{thm}

By the definition of K-stability for pairs which will be recalled in the next section, the Futaki invariant depends linearly on the angle $\beta$. Thus Theorem \ref{thm} leads immediately to  the following

\begin{cor}\label{cor}
If $X$ is K-stable(semi-stable), then for any smooth anti-canonical divisor $D$, the pair $(X, D)$ is K-stable(semi-stable) with respect to angle $\beta\in(0,1]$.
\end{cor}

This corollary provides evidence to the picture described in \cite{Do3} that a smooth K\"ahler-Einstein metric on $X$ should come from a complete Calabi-Yau metric on $X\setminus D$ by increasing the angle from $0$ to $2\pi$. The relevant definitions will be given in the next section. The strategy to prove K-unstability for negative angles is by studying a particular test configuration, namely the deformation to the normal cone of $D$. To deal with the zero angle case we shall construct ``approximately balanced" embeddings using the Calabi-Yau metric on $D$.  In \cite{Od}, Odaka proved that a Calabi-Yau manifold is K-stable, by a purely algebro-geometric approach. It is very likely that his method can give an alternative proof of the above theorem, but the one we take seems to be more quantitative.\\

\section{K-stability for pairs}
We first recall the definition of K-stability. 

\begin{defi} Let $(X,L)$ be a polarized manifold. A \emph{test configuration for $(X, L)$}  is a $\C^*$ equivariant flat family  $(\X, \L)\rightarrow\C$
such that $(\X_1, \L_1)$ is isomorphic to $(X, L)$. $(\X,  \L)$ is called \emph{trivial} if it is isomorphic to the product $(X, L)\times \C$ with the trivial action on $(X, L)$ and the standard action on $\C$.
\end{defi}

Suppose $D$ is a smooth divisor in $X$, then any test configuration $(\X, \L)$ induces a  test configuration  $(\D, \L)$ by simply taking the flat limit of the $\C^*$ orbit of $D$ in $\X_1$. We call $(\X, \D, \L)$ a \emph{test configuration for $(X, D, L)$}.   Given any test configuration $(\X, \D, \L)$ for $(X, D, L)$, we denote by $A_k$ and $\tilde{A}_k$ the infinitesimal generators for the $\C^*$ action on $H^0(\X_0, \L^k_0)$ and $H^0(\D_0, \L_0^k)$ respectively.  By general theory for $k$ large enough we have the following expansions
$$d_k:=h^0(\X_0, \L_0^k)=a_0 k^n+a_1 k^{n-1}+O(k^{n-2}), $$
$$w_k:=tr(A_k)=b_0 k^{n+1}+b_1k^n+O(k^{n-1}), $$
$$\tilde{d}_k:=h^0(\D_0, \L_0^k)=\tilde{a}_0 k^{n-1}+\tilde{a}_1 k^{n-2}+O(k^{n-3}), $$
$$\tilde{w}_k:=tr(\tilde{A}_k)=\tilde{b}_0 k^{n}+\tilde{b}_1k^{n-1}+O(k^{n-2}).$$

\begin{defi}
For any real number $\beta$, the \emph{Futaki invariant} of a test configuration $(\X, \D, \L)$ with respect to angle $\beta$ is 
$$Fut(\X, \D, \L, \beta)=\frac{2(a_1b_0-a_0b_1)}{a_0}+(1-\beta)(\tilde{b}_0-\frac{\tilde{a}_0}{a_0}b_0).$$ 
\end{defi}

When $\beta=1$  we get the usual Futaki invariant of a test configuration  $(\X, \L)$ 
$$Fut(\X, \L)=\frac{2(a_1b_0-a_0b_1)}{a_0}.$$

\begin{defi}
A polarized manifold  $(X, L)$ is called \emph{K-stable(semistable)} if $Fut(\X, \L)>0$$(\geq 0)$ for any nontrivial test configuration $(\X, \L)$. Similarly,  $(X, D, L)$ is called \emph{K-stable(semistable) with respect to angle $\beta$} if $Fut(\X, \D, \L, \beta)>0$$(\geq 0)$ for any nontrivial test configuration $(\X, \D, \L)$.
\end{defi}

When the central fiber $(\X_0, \D_0)$ is smooth,  by Riemann-Roch
the Futaki invariant then has a differential geometric  expression as
$$Fut(\X, \D, \L, \beta)=\int_{\X_0} (S-\underline{S})H \frac{\omega^n}{n!}-(1-\beta)(\int_{\D_0} H\frac{\omega^{n-1}}{(n-1)!}-\frac{Vol(\D_0)}{Vol(\X_0)}\int_{\X_0} H\frac{\omega^{n}}{n!}), $$
where $\omega$ is an $S^1$ invariant K\"ahler metric in $2\pi c_1(\L_0)$ and $H$ is the Hamiltonian function generating the $S^1$ action on $\L_0$. This differs from the usual Futaki invariant by an extra term which reflects the cone angle.\\

The above abstract notion of K-stability  is closely related to Chow stability for projective varieties, which we now recall. 
Given a $\C^*$ action on $\C\P^N$, and suppose the induced $S^1$ action preserves the Fubini-Study metric. Then the infinitesimal generator is given by a Hermitian matrix, say $A$. The Hamiltonian function for the $S^1$ action on $\C\P^N$ is 
$$H_A(z)=\frac{z^*Az}{|z|^2}.$$
Given a projective manifold $V$ in $\C\P^N$, we define the \emph{center of mass} of $V$
$$\mu(V)=\int_V\frac{zz^*}{|z|^2}d\mu_{FS}-\frac{Vol(V)}{N+1}Id\in \sqrt{-1} \mathfrak{su}(N+1), $$  viewing $\C\P^N$ as a co-adjoint orbit in $\mathfrak{su}(N+1)$. 
Define the \emph{Chow weight of $V$ with respect to $A$}  to be
$$CH(V, A)=-Tr(\mu(V)\cdot A)=-\int_V H_A d\mu_{FS}+\frac{Vol(V)}{N+1}TrA.$$
Notice this vanishes if $A$ is a scalar matrix. The definition is not sensitive to singularities of $V$ so one may define the Chow weight of any algebraic cycles in a natural way. It is well-known that the $CH(e^{tA}.V, A)$  is a decreasing function of $t$, see for example \cite{Do2}. So
 \begin{equation}\label{monotonicity}
 CH(V, A)\leq CH(V_{\infty}, A), 
 \end{equation}
 where $V_{\infty}$ is the limiting Chow cycle of $e^{tA}.V$ as $t\rightarrow-\infty$. $V_{\infty}$ is fixed by the $\C^*$ action and then $CH(V_{\infty}, A)$ is an algebraic geometric notion, i.e. independent of the Hermitian metric we choose on $\C^{N+1}$.\\

This well-known theory readily extends to  pairs, see \cite{Do3}, \cite{AD}. We consider a pair of varieties $(V, W)$ in $\C\P^N$ where $W$ is a subvariety of $V$. Given a parameter $\lambda\in[0,1]$, we define the \emph{ center of mass of $(V, W)$ with parameter $\lambda$} 
$$\mu(V, W, \lambda)=\lambda\int_V \frac{zz^*}{|z|^2}d\mu_{FS}+(1-\lambda)\int_W \frac{zz^*}{|z|^2} d\mu_{FS}-\frac{\lambda Vol(V)+(1-\lambda)Vol(W)}{N+1}Id, $$ and the \emph{
Chow weight}  with parameter $\lambda$:
$$CH(V, W, A, \lambda)=-Tr(\mu(V, W, \lambda)\cdot A).$$
A pair $(V, W)$ with vanishing center of mass with parameter $\lambda$ is called a \emph{$\lambda$-balanced embedding}. \\

Now given a test configuration $(\X, \D, \L)$, it is explained in \cite{RT1} and \cite{Do2}(see also \cite{PS}) that for $k$ large enough one can realize it by a family of projective schemes in $\P(H^0(X, L^k)^*)$ with a one parameter group action. Moreover one could arrange that the fiber $(\X_1, \D_1, \L_1)$ is embedded into $\P(H^0(X, L^k)^*)$ with a prescribed Hermitian metric, and the $\C^*$ action  is generated by a Hermitian matrix $-A_k$(negative sign because we are taking the dual). Then as in \cite{Do2}  the Futaki invariant is the limit of Chow weight:

\begin{equation} \label{ChowFutaki}
\lim_{k\rightarrow\infty}k^{-n}CH_k(\X_0, \D_0, -A_k, \lambda))=Fut(\X, \D, \L, \beta), 
\end{equation} with $\beta=\frac{3\lambda-2}{\lambda}$.\\

\section{Proof of the main theorem}
From now  on we assume $X$ is a Fano manifold of dimension $n$, $D$ is a smooth anti-canonical divisor and  the polarization is given by $L=-K_X$.
We first prove the part of unstability in theorem \ref{thm},  by considering the deformation to the normal cone of $D$, as studied by Ross-Thomas \cite{RT2}. We blow up $D\times \{0\}$ in the total space $X\times\C$ and get a family $\pi:\X\rightarrow\C$. The exceptional divisor $P$ is equal to the projective completion $\P(\nu_D\oplus \C)$ of the normal bundle $\nu_D$ in $X$. The central fiber $\X_0$ is the gluing of $P$ to $X$ along $D=\P(\nu_D)$.  There is a $\C^*$ action on $\X$ coming from the trivial action on $X$ and the standard $\C^*$ action on $\C$. Let $\D$ be the proper transform of $D\times\C$. This is $\C^*$ invariant, and its intersection with  the central fiber is the zero section $\P(\C)\subset \P(\nu_D\oplus\C)$(The readers are referred to \cite{RT2} for a very nice picture of a deformation to the normal cone). The line bundle we use is $\L_c=L(-cP)$($c$ is rational).   It is shown in \cite{RT2} that $\L_c$ is ample when $c\in (0,1)$. There is also a natural lift of the $\C^*$ action to $\L_c$, so that we get test configurations $(\X, \D, \L_c)$ parametrized by $c$.  We follow \cite{RT2}  to compute the Futaki invariant.  Pick a sufficiently large integer $k$ so that $ck$ is an integer. 
We have the decomposition
\begin{equation*}H^0(\X, \L_c^k)=\bigoplus_{i=1}^{ck}t^{ck-i}H^0(X, L^{k-i})\oplus t^{ck}\C[t]H^0(X, L^k), \end{equation*}
where $t$ is the standard holomorphic coordinate on $\C$. 
Using the short exact sequence $$0\rightarrow H^0(X, L^{i-1})\rightarrow H^0(X, L^{i})\rightarrow H^0(D, L^i)\rightarrow0, $$
we obtain
\begin{eqnarray*}
H^0(\X_0, \L_c^k)&=&H^0(\X, \L_c^k)/tH^0(\X, \L_c^k)\\&=&H^0(X, L^{(1-c)k} )\oplus \bigoplus_{i=0}^{ck-1} t^{ck-i} H^0(D, L^{k-i}).
\end{eqnarray*}
This is indeed the weight decomposition of $H^0(\X_0, \L_c^k)$ under the $\C^*$ action. Note the weight is $-1$ on $t$. So 
$$\dim H^0(\X_0, \L_c^k)=\dim H^0(X, L^{(1-c)k})+\sum_{i=0}^{ck-1} \dim H^0(D, L^{k-i})=\dim H^0(X, L^{ck}).$$
This actually shows the flatness of the family $(\X, \D, \L)$. Thus by Riemann-Roch, 
$$a_0=\frac{1}{n!}\int_X c_1(L)^n, $$ and 
 $$a_1=\frac{1}{2(n-1)!}\int_Xc_1(-K_X)\cdot c_1(L)^{n-1}=\frac{na_0}{2}.$$
 The weight is given by 
\begin{eqnarray*}
w_k&=& -\sum_{i=0}^{ck-1} (ck-i)\dim H^0(D, L^{k-i})\\
       &=& -\sum_{i=0}^{ck-1} (ck-i)(\frac{(k-i)^{n-1}}{(n-1)!} \int_D c_1(L)^{n-1}+O(k^{n-3}))\\&=&
       -na_0\int_0^c (c-x)(1-x)^{n-1}dx\cdot k^{n+1}-\frac{nca_0}{2}k^n+O(k^{n-1}).
\end{eqnarray*}
 So $$b_0=(\frac{1-(1-c)^{n+1}}{n+1}-c)a_0, $$
 and $$b_1=-\frac{nca_0}{2}.$$
 Thus the ordinary Futaki invariant for the test configuration $(\X, \L)$ is given by 
 $$Fut_c(\X, \L)=\frac{2(a_1b_0-a_0b_1)}{a_0}=n(\frac{1-(1-c)^{n+1}}{n+1})a_0.$$
 Note
 $$H^0(\D, \L_c^k)=H^0(D\times\C, L^k\otimes (t)^{ck})=t^{ck} \C[t] H^0(D, L^k).$$
 So  $$H^0(\D_0, \L_c^k)=H^0(\D, \L_c^k)/tH^0(\D, \L_c^k)=t^{ck}H^0(D, L^k).$$
 Thus we see
 $$\tilde{a}_0=\int_D \frac{c_1(L)^{n-1}}{(n-1)!}=na_0, $$
 and $$\tilde{b}_0=-c\int_D\frac{c_1(L)^{n-1}}{(n-1)!}=-nca_0.$$
Therefore,
 \begin{eqnarray*}
 Fut_c(\X, \D, \L, \beta)&=&Fut_c(\X,\L)+(1-\beta)(\tilde{b}_0-\frac{\tilde{a}_0}{a_0}{b_0})\\
 &=& [n(\frac{1-(1-c)^{n+1}}{n+1})+(1-\beta)(-nc+n(c-\frac{1-(1-c)^{n+1}}{n+1}))]a_0\\
 &=&n\beta\frac{1-(1-c)^{n+1}}{n+1}a_0.
 \end{eqnarray*}
 Therefore  for $\beta<0$  this particular test configuration gives rise to unstability, and for $\beta=0$  the pair $(X, D)$ can not be stable.\\
 
Now we move on to prove K-semistability for $\beta=0$.
 Using again the short exact sequence 
$$0\rightarrow H^0(X, L^{j-1})\rightarrow H^0(X, L^{j})\rightarrow H^0(D, L^j)\rightarrow 0$$ successively  we can choose a splitting
 \begin{equation} \label{splitting}H^0(X, L^k)=H^0(X, L^{s-1})\oplus \bigoplus_{j=s}^{k}H^0(D, L^j)\end{equation} for $s$ large enough and all $k>s$. By Yau's  theorem \cite{Ya1} there is a unique Ricci flat metric $\omega_0$ in $c_1(L)|_D$. This defines a Hermitian metric on $H^0(D, L^j)$ by the $L^2$ inner product.  We can put an arbitrary metric on $H^0(X, L^{s-1})$, and make the splitting (\ref{splitting}) orthogonal. We also identify the vector spaces with their duals using these metrics. Take $s$ large enough so that $D$ embeds into $\P(H^0(D, L^j))$ and $X$ embeds into $\P(H^0(X, L^j))$ for all $j\geq s-1$. Choosing an orthonormal basis of $H^0(D, L^j)$ we get an embedding $f_j: D\rightarrow \P(H^0(D, L^j))\cong \P^{n_j-1}$(here $n_j=\dim H^0(D, L^j)$). We also pick an arbitrary embedding $f_{s-1}: X\rightarrow\P(H^0(X, L^{s-1}))$. Denote by $D_j$ the image of $f_j$, and let $N(D_{j-1}, D_j)$ be the variety consisting of all points  in $\P(H^0(D, L^{j-1})\oplus H^0(D, L^j))\subset \P(H^0(X, L^k))$ of the form $[u f_{j-1}(p): v f_j(p)]$ for $p\in D$ and  $u, v\in\C.$ The projection map $\pi_j: N(D_{j-1}, D_j)\rightarrow D$ makes it a $\P^1$ bundle over $D$. This is isomorphic to the projective completion of the normal bundle of $D$ in $X$. Let $X_k$ be  the union of all these $N(D_{j-1}, D_j)$($s\leq j\leq k$) together with $f_{s-1}(X)$. Then it is not hard to see that as a pair of Chow cycles $(X_k, D_k)$ lies in the closure of the $PGL(d_k;\C)$ orbit of a smooth embedding of $(X, D)$ in $\P(H^0(X, L^k))$.  We want to estimate its center of mass. The following two lemmas involve some calculation and the proof will be deferred to the end of this section.
 
 \begin{lem} \label{lemma1} For $s\leq j\leq k$ we have
 $${\pi_j}_*\omega_{FS}^{n}=\sum_{i=0}^{n-1}\omega_{j}^{i}\wedge\omega_{j-1}^{n-1-i}, $$
where $\omega_j=f_j^*\omega_{FS}$.
 \end{lem} 
 This lemma implies  that $$Vol(N(D_{j-1}, D_j))=\frac{1}{n!}\sum_{i=0}^{n-1}j^i(j-1)^{n-1-i}\cdot (n-1)!Vol(D)=(j^n-(j-1)^n)Vol(X). $$ Summing over $j$ we see that $Vol(X_k)=k^n Vol(X)$.\\
 
 Notice $N(D_{j-1}, D_j)$ can only contribute to the $H^0(D, L^{j-1})$ and $H^0(D, L^j)$ components of the center of mass of $X_k$. Denote by $Z_j=(Z_j^1, \cdots, Z_j^{n_j})$ the homogeneous coordinates on $H^0(D, L^j)$ for $s\leq j\leq k$, and by $Z_{s-1}$ the homogeneous coordinate on $H^0(X, L^{s-1})$. Then we have 

 \begin{lem} \label{lemma2} For $s\leq j\leq k$ we have
$$ {\pi_j}_* \frac{Z_jZ_{j-1}^*}{|Z_j|^2+|Z_{j-1}|^2}\omega_{FS}^n=0, $$
$$ {\pi_j}_* \frac{Z_{j-1}Z_{j}^*}{|Z_j|^2+|Z_{j-1}|^2}\omega_{FS}^n=0, $$
$$ {\pi_j}_* \frac{Z_jZ_{j}^*}{|Z_j|^2+|Z_{j-1}|^2}\omega_{FS}^n=\frac{Z_jZ_j^*}{|Z_j|^2}\cdot\sum_{i=0}^{n-1}\frac{i+1}{n+1}\omega_j^i\wedge \omega_{j-1}^{n-1-i}, $$
$$ {\pi_j}_* \frac{Z_{j-1}Z_{j-1}^*}{|Z_j|^2+|Z_{j-1}|^2}\omega_{FS}^n=\frac{Z_{j-1}Z_{j-1}^*}{|Z_{j-1}|^2}\cdot\sum_{i=0}^{n-1}\frac{n-i}{n+1}\omega_j^i\wedge \omega_{j-1}^{n-1-i}, $$
\end{lem} 

 This lemma implies that the center of mass $\mu(X_k)$ also splits as the direct sum of $\mu_j$'s. For $j$ between $s$ and $k-1$ we have
 $$\mu_j(X_k)=\int_{X_k} \frac{Z_j Z_j^*}{|Z|^2} \frac{\omega_{FS}^n}{n!}=\frac{1}{n!}\int_D \frac{Z_jZ_j^*}{|Z_j|^2} \sum_{i=0}^{n-1}(\frac{i+1}{n+1}\omega_j^i\wedge\omega_{j-1}^{n-1-i}+
 \frac{n-i}{n+1}\omega_{j+1}^i\wedge\omega_{j}^{n-1-i}), $$
 while $$\mu_k(X_k)=\int_{X_k} \frac{Z_k Z_k^*}{|Z|} \frac{1}{n!}\omega_{FS}^n=\frac{1}{n!}\int_D \frac{Z_kZ_k^*}{|Z_k|^2} \sum_{i=0}^{n-1}\frac{i+1}{n+1}\omega_k^i\wedge\omega_{k-1}
 ^{n-1-i}, $$
and $$\mu_{s-1}(X_k)=\frac{1}{n!}\int_D \frac{Z_{s-1}Z_{s-1}^*}{|Z_{s-1}|^2} \sum_{i=0}^{n-1}\frac{n-i}{n+1}\omega_s^i\wedge \omega_{s-1}^{n-1-i}+\int_{X_{s-1}} \frac{Z_{s-1}Z_{s-1}^*}{|Z_{s-1}|^2
}\frac{\omega_{FS}^n}{n!}.$$ \\
The induced metric $\omega_j$ is related to the original metric $\omega_0$  by the ``density of state" function:
$$\omega_j=j\omega_0+\sqrt{-1}\p\bp\log \rho_j(\omega_0). $$
It is well-known that we have the following expansion(see \cite{Ca}, \cite{Ze}, \cite{Lu}, \cite{MM}) $$\rho_j(\omega_0)=j^{n-1}+\frac{S(\omega_0)}{2}j^{n-2}+O(j^{n-3})=j^{n-1}+O(j^{n-3}), $$ since $\omega_0$ is Ricci flat.  Thus
$$\omega_j^i\omega_{j-1}^{n-1-i}=j^i(j-1)^{n-1-i}\omega_0^{n-1}(1+O(j^{-3})).$$
To estimate $\mu_j$ recall we have chosen an orthonormal basis $\{s_j^l\}$ of $H^0(D, L^j)$ and  we can assume  $\mu_j$ is a diagonal matrix.
 Then for $s\leq j\leq k-1$ we obtain
 $$\mu_j^l(X_k)=\int_D \frac{|s_j^l|^2(1+O(j^{-3}))}{j^{n-1}+O(j^{n-3})}\sum_{i=0}^{n-1}(\frac{i+1}{n+1}j^i(j-1)^{n-1-i}+\frac{n-i}{n+1}(j+1)^ij^{n-1-i})\frac{\omega_0^{n-1}}{n!}.$$
 It is easy to see that
 $$\sum_{i=0}^{n-1}(\frac{i+1}{n+1}j^i(j-1)^{n-1-i}+\frac{n-i}{n+1}(j+1)^ij^{n-1-i})=nj^{n-1}+O(j^{n-3}).$$
 Thus 
 $$\mu_j^l(X_k)=1+O(j^{-2}).$$
 For $j=k$, we have $$\mu_k^l(X_k)=1/2+O(k^{-1}).$$
 For $j=s-1$, we have $$\mu_{s-1}^l(X_k)=O(1).$$
 The center of mass of the pair $(X_k, D_k)$ with respect to $\lambda=2/3$ is given by
 $$\mu(X_k, D_k, 2/3)=\frac{2}{3}\mu(X_k)+\frac{1}{3}\mu(D_k)-\underline{\mu} \cdot Id,$$
 where we denote  $$\underline{\mu}=\frac{2Vol(X_k)+Vol(D_k)}{3d_k}=\frac 23+O(k^{-2}).$$ 
 Thus for $s\leq j\leq k-1$ and $0\leq l\leq n_j$ we have 
 $$\mu_j^l(X_k, D_k, 2/3)=O(j^{-2})+O(k^{-2}).$$
 Since $n_j$ is a polynomial  of degree $n-1$ in $j$, we obtain
 $$|\mu_j(X_k, D_k, 2/3)|_2=(\sum_{l=0}^{n_j} |\mu_j^l(X_k, D_k, 2/3)|^2)^{1/2}=O(j^{\frac{n-5}{2}}), $$
 and $$\sum_{j=s}^{k-1} |\mu_j(X_k, D_k, 2/3)|_2=O(k^{\frac{n-3}{2}}).$$
 For $j=k$, 
 we have 
 $$\mu_k^l(D_k)=\int_D\frac{|s_k^l|^2}{k^{n-1}+O(k^{n-3})}(1+O(k^{-2}))\frac{k^{n-1}\omega_0^{n-1}}{(n-1)!}=1+O(k^{-2}).$$
So $$\mu_k^l(X_k, D_k)=O(k^{-1}), $$
and $$|\mu_k(X_k, D_k)|_2=O(k^{\frac{n-3}{2}}).$$
Therefore we obtain $$|\mu(X_k, D_k)|_2=O(k^{\frac{n-3}{2}}).$$
So for a smoothly embedded $(X, D)$ in $\P(H^0(X, L^k))$ we have $$\inf_{g\in PGL(d_k;\C)} |\mu(g.(X, D))|_2=O(k^{\frac{n-3}{2}}).$$
In particular there are embeddings $\iota_k: (X, D)\rightarrow \P(H^0(X, L^k))$ such that $$|\mu(\iota_k(X, D))|_2=O(k^{\frac{n-3}{2}}).$$
 Now any test configuration $(\X, \D, \L)$ can be represented by a family in $\P(H^0(X, L^k))$ such that the fiber $(\X_1, \D_1, \L_1)$ is embedded by $\iota_k$ and the $\C^*$ action is generated by a Hermitian matrix $A_k$. Again by general theory $|A_k|_{2}^2=TrA_k^2=O(k^{n+2})$.   Therefore  by monotonicity of the Chow weight we obtain
 \begin{eqnarray*}
 CH_k(\X_0, \D_0, -A_k, 2/3)&\geq& CH_k(\X_1, \D_1, -A_k, 2/3)\\&\geq& -\inf_{g\in PGL(d_k;\C)}|\mu(g.(X, D))|_{2}\cdot|-A_k|_{2}\\&\geq &-O(k^{n-\frac{1}{2}}).
 \end{eqnarray*}
 Thus by (\ref{ChowFutaki})  $$Fut(\X, \D, \L, 0)=\lim_{k\rightarrow\infty} k^{-n}CH_k(\X_0, \D_0, -A_k, \frac{2}{3})\geq 0.$$ This finishes the proof of Theorem \ref{thm}.\\

Now we prove Lemmas \ref{lemma1} and \ref{lemma2}.
 In general suppose  there are two embeddings $f_1:D\rightarrow \P^l$ and $f_2: D\rightarrow \P^m$.  As before, let $N(D)$ be the variety in  $\P^{l+m+1}$ containing all points of the form $(tf_1(x), sf_2(x))$ where $t, s\in \C$. Intuitively $N(D)$ is ruled by all lines connecting $f_1(x)$ and $f_2(x)$ for $x\in D$.  Choose a local coordinate chart $U$ in $D$ such that the image $f_1(U)$ and $f_2(U)$ are contained in a standard coordinate chart for the projective spaces $\P^l$ and $\P^m$ respectively. Let $[1:z]$ and $[1:w]$ be local coordinates in $\P^l$ and $\P^m$. 
 Under unitary transformations we may assume $f_1(x_0)=[1:0]$ and $f_2(x_0)=[1:0]$. The line connecting $f_1(x_0)$ and $f_2(x_0)$ is parametrized as $[1: 0:t: 0]$ for $t\in \C$. Along this line we have
\begin{eqnarray*}
\omega_{FS}&=& \frac{\sqrt{-1}}{2\pi}\p\bp \log(1+|z|^2+|t|^2+|t|^2|w|^2)\\                   
                       &=&\frac{\sqrt{-1}}{2\pi}\cdot\frac{(1+|t|^2)\sum_idz^i\wedge d\bz^i+|t|^2(1+|t|^2)\sum_jdw^j\wedge d\bw^j+dt\wedge d\bt}{(1+|t|^2)^2}.
\end{eqnarray*}
Thus 
$$\omega_{FS}^n=n(\frac{\sqrt{-1}}{2\pi})^n (1+|t|^2)^{-n-1}(\sum_idz^i\wedge d \bz^i+|t|^2\sum_jdw^j\wedge d\bw^j)^{n-1}\wedge dt\wedge d\bt.$$ Hence integrating along the $\P^1$ we get 
\begin{eqnarray*}
\int_{\P^1} \omega_{FS}^n&=&\frac{1}{2\pi} \int_{\C} n (\omega_1+|t|^2\omega_2)^{n-1}\wedge (1+|t|^2)^{-n-1}\sqrt{-1}dt\wedge d\bt\\
&=& \frac{1}{2\pi}\int_0^{\infty} n \sum_{j=0}^{n-1}\binom{n-1}{j}\omega_1^j\wedge \omega_2^{n-1-j} x^j(1+x)^{-n-1}dx\\
&=&\sum_{j=0}^{n-1}\omega_1^j\wedge \omega_2^{n-1-j}.
\end{eqnarray*}
This proves lemma \ref{lemma1}.\\

For the center of mass we compute
$$\int_{\P^1} \frac{1}{1+|t|^2}\omega_{FS}^n=\sum_{j=0}^{n-1}\frac{j+1}{n+1}\omega_1^j\wedge \omega_2^{n-1-j}, $$
and $$\int_{\P^1}\frac{|t|^2}{1+|t|^2}\omega_{FS}^n=\sum_{j=0}^{n-1}\frac{n-j}{n+1}\omega_1^j\wedge \omega_2^{n-1-j}.$$
Thus globally we obtain
$$\int_{N(D)}\frac{zz^*}{|z|^2+|w|^2} \omega_{FS}^n=\int_{D} \frac{zz^*}{|z|^2} \sum_{j=0}^{n-1}\frac{j+1}{n+1}\omega_1^j\wedge \omega_2^{n-1-j}, $$
and 
$$\int_{N(D)}\frac{ww^*}{|z|^2+|w|^2}\omega_{FS}^n=\int_D \frac{ww^*}{|w|^2}\sum_{j=0}^{n-1}\frac{n-j}{n+1}\omega_1^j\wedge \omega_2^{n-1-j}.$$
Also notice by symmetry of $N(D)$ under the map $w\mapsto-w$ we have
$$\int_{N(D)}\frac{zw^*}{|z|^2+|w|^2}\omega_{FS}^n=0.$$ Similarly
$$\int_{N(D)}\frac{wz^*}{|z|^2+|w|^2}\omega_{FS}^n=0.$$ 
This proves lemma \ref{lemma2}.\\

\begin{rmk}
In the case when $X$ is $\P^1$ and $D$ consists of two points,  one can indeed find the precise balanced embedding for $\lambda=2/3$.
 In $\P^k$ let $L$ be the chain of lines $L_i$ connecting $p_i$ and $p_{i+1}$($0\leq i\leq k-1$), where $p_i$ is the $i$-th coordinate point. Then it is easy to see that $L$ is the degeneration limit of a smooth degree $k$ rational curve, and it is exactly $\frac{2}{3}$-balanced. It is well-known that a rational normal curve in  $\P^k$ is always Chow polystable, it follows by linearity that it is also Chow polystable for $\lambda\in (2/3, 1]$.
\end{rmk}

\noindent {\bf Acknowledgements:} The author is grateful to Professor Simon Donaldson for enlightening discussions. He would also like to thank Professor Xiuxiong Chen for constant encouragement. This work was partly supported  by  a postdoc grant under European Research Council award No 247331.


\begin{thebibliography}{0}
\bibitem{AD} C. Arezzo, A. Della Vedova, Notes in preparation. 
\bibitem{Ca} D. Catlin, \emph{The Bergman kernel and a theorem of Tian}, in ÔAnalysis and geometry in several complex variablesÕ (Katata, 1997) 1-23, Birhauser, Boston, 1999.
\bibitem{Do1} S. K. Donaldson. \emph{Scalar curvature and stability of toric varieties}. J. Differential. Geom. 62 (2002), no. 2, 289-349.
\bibitem{Do2} S. K. Donaldson. \emph{Lower bounds on the Calabi functional}, J. Differential Geom. 70 (2005), no. 3, 453--472.
\bibitem{Do3} S. K. Donaldson. \emph{Discussion of the K\"ahler-Einstein problem}, notes available on http://www2.imperial.ac.uk/~skdona/.
\bibitem{Do4} S. K. Donaldson. \emph{K\"ahler metrics with cone singularities along a divisor,}  arXiv:1102.1196. 
\bibitem{Li} C. Li. \emph{Remarks on logarithmic K-stability}, arxiv: 1104.0428. 
\bibitem{Lu} Z-Q. Lu. \emph{On the lower order terms of the asymptotic expansion of Tian-Yau-Zelditch}, Amer. J. Math, 122(2000), 235--273.
\bibitem{Ma} T. Mabuchi. \emph{K-stability of constant scalar curvature polarization}, arxiv: 0812.4093.
\bibitem{MM} M. Marinescu, X-N. Ma. \emph{Holomorphic Morse inequalities and Bergman kernels},  Progress in Mathematics, 254. Birkh\"auser Verlag, Basel, 2007.
\bibitem{Od} Y. Odaka. \emph{The Calabi conjecture and K-stability
}, arxiv: 1010.3597.
\bibitem{PS} D. Phong, J. Sturm. \emph{Test configurations for K-stability and geodesic rays.} J. Symplectic Geom. 5 (2007), no. 2, 221--247.
\bibitem{RT1} J. Ross, R. Thomas. \emph{A study of Hilbert-Mumford criterion for the stability of projective varieties}, J. Algebraic. Geom. 16(2007), 201--255.
\bibitem{RT2} J. Ross, R. Thomas. \emph{An obstruction to the existence of constant scalar curvature K\"ahler metrics}, J. Differential. Geom. 72(2006), 429--466.
\bibitem{St} J. Stoppa. \emph{K-stability of constant scalar curvature K\"ahler manifolds}. Adv. Math. 221 (2009), no. 4, 1397--1408.
\bibitem{Sz} G. Sz\'ekelyhidi, \emph{Extremal metrics and K-stability}, Ph.D. thesis, 2006.
\bibitem{Ti} G. Tian. \emph{K\"ahler-Einstein metrics with positive scalar curvature}, Invent. Math. 130 (1997), 1--39.
\bibitem{TY} G. Tian, S-T. Yau. \emph{Complete K\"ahler manifolds with zero Ricci curvature}. I. J. Amer. Math. Soc. 3 (1990), no. 3, 579--609.
\bibitem{Ya1} S-T. Yau. \emph{On the Ricci curvature of a compact K\"ahler manifold and the complex Monge-Amp\`ere equation}, I. Comm. Pure. Appl. Math. 31(1978), no. 3, 339--411.
\bibitem{Ya2} S-T. Yau. \emph{Open problems in geometry}, Proc. Sympos. Pure Math. 54 (1993),
1--28.
\bibitem{Ze} S. Zelditch. \emph{Szeg\"o kernel and a theorem of Tian},
Int. Math. Res. Notices 6(1998), 317--331.
\end{thebibliography}
\end{document}